\documentclass[journal]{IEEEtran}
\usepackage{amsfonts}
\usepackage{amsmath}
\usepackage{amssymb}
\usepackage{epsfig}
\usepackage{graphics}
\usepackage{graphicx}

\epsfclipon
\newtheorem{theorem}{Theorem}[section]

\newtheorem{corollary}[theorem]{Corollary}

\newtheorem{definition}{Definition}[section]
\newtheorem{example}{Example}[section]

\newtheorem{lemma}[theorem]{Lemma}

\newtheorem{problem}{Problem}

\newtheorem{remark}{Remark}[section]

\input{tcilatex}
\begin{document}

\title{Structural Analysis of Viral Spreading Processes in Social and
Communication Networks Using Egonets}
\author{Victor M. Preciado,~\IEEEmembership{Member,~IEEE,} Moez Draief, and
Ali Jadbabaie,~\IEEEmembership{Senior~Member,~IEEE} \thanks{%
Manuscript Received ---------.} \thanks{%
V.M. Preciado and A. Jadbabaie are with the Department of Electrical and
Systems Engineering at the University of Pennsylvania, Philadelphia, PA
19104 USA. (e-mail: preciado@seas.upenn.edu; jadbabai@seas.upenn.edu).} 
\thanks{%
M. Draief is with the Department of Electrical and Electronic Engineering at
Imperial College, London, SW7 2AZ UK. (e-mail: m.draief@imperial.ac.uk). } 
\thanks{%
This work was supported by ONR MURI \textquotedblleft Next Generation
Network Science\textquotedblright\ and AFOSR \textquotedblleft Topological
And Geometric Tools For Analysis Of Complex Networks\textquotedblright .} 
\thanks{%
MD is supported by QNRF through grant NPRP 09-1150-2-448. MD holds a
Leverhulme Trust Research Fellowship RF/9/RFG/2010/02/08.}}
\maketitle

\begin{abstract}
We study how the behavior of viral spreading processes is influenced by
local structural properties of the network over which they propagate. For a
wide variety of spreading processes, the largest eigenvalue of the adjacency
matrix of the network plays a key role on their global dynamical behavior.
For many real-world large-scale networks, it is unfeasible to exactly
retrieve the complete network structure to compute its largest eigenvalue.
Instead, one usually have access to myopic, egocentric views of the network
structure, also called \emph{egonets}. In this paper, we propose a
mathematical framework, based on algebraic graph theory and convex
optimization, to study how local structural properties of the network
constrain the interval of possible values in which the largest eigenvalue
must lie. Based on this framework, we present a computationally efficient
approach to find this interval from a collection of egonets. Our numerical
simulations show that, for several social and communication networks, local
structural properties of the network strongly constrain the location of the
largest eigenvalue and the resulting spreading dynamics. From a practical
point of view, our results can be used to dictate immunization strategies to
tame the spreading of a virus, or to design network topologies that
facilitate the spreading of information virally.
\end{abstract}


\section{Introduction}

Understanding the behavior of viral spreading processes taking place in
large complex networks\ is of critical interest in mathematical epidemiology 
\cite{New03,BLMCH06}. Spreading processes are relevant in many real
scenarios, such as disease spreading in human populations \cite{AB00}--\cite%
{Het00}, malware propagation in computer networks \cite{BBCS05}--\cite{Kle07}%
, or information dissemination in online social networks \cite{GGLT04}--\cite%
{LAH06}. To study viral spreading processes, a variety of stochastic
dynamical models has been proposed in the literature \cite{GMT05}--\cite%
{PCFVF11}. In these models, the steady-state infection of the network
presents two different regimes depending on the virulence of the infection
and the structure of the network of contacts. In one of the regimes, an
initial infection dies out at a fast (usually exponential) rate. In the
other regime, an initial infection becomes an epidemic. Both numerical and
analytical results show that these two regimes are separated by a phase
transition at an epidemic threshold determined by both the virulence of the
infection and the topology of the network. One of the most fundamental
questions in mathematical epidemiology is to find the value of the epidemic
threshold in terms of the virus model and the contact network.

In many cases of practical interest it is unfeasible to exactly retrieve the
complete structure of a network of contacts. In these cases, it is
impossible to exactly compute the epidemic threshold. On the other hand, in
most cases one can easily retrieve the structure of egocentric views of the
network, also called \emph{egonets}\footnote{%
A rigorous definition of egonet, in graph-theoretical terms, will be given
in Section \ref{Spectral Analysis from Egonets}.}. To estimate the value of
the epidemic threshold, researchers have proposed a variety of random
network models in which they can prescribe structural properties that can be
retrieved from these egonets, such as the degree distribution \cite%
{NSW01,CLV03}, local correlations \cite{PVV01,New02b}, or clustering \cite%
{New09}.

Although random networks are the primary tool to study the impact of local
structural features on the epidemic threshold \cite{DM10}, this approach
presents a major flaw: Random network models implicitly induce many
structural properties that are not directly controlled but can have a strong
influence on the value of the epidemic threshold. For example, it is
possible to find two networks having the same degree distribution, but with
opposite dynamical behavior\emph{\ }\cite{ALWD05}. Therefore, it is
difficult (if not impossible) to isolate the role of a particular structural
property in the network performance using random network models.
Furthermore, many real networks present weighted edges representing, for
example, bandwidth in communication networks or resistance in electric
networks. Current random networks fail to faithfully recover both the
structure of the network and the distribution of weights over the links. In
this paper, we develop a mathematical framework, based on algebraic graph
theory and convex optimization, to study how the structure of local egonets
constrain the interval of possible values in which the epidemic threshold
must lie. As a result of our analysis, we present a computationally
efficient approach to find this interval from a collection of egonets
extracted from a (possibly) weighted network. Our numerical simulations show
that the resulting interval is very narrow for several social and
communication networks. This illustrates the fact that, for many real
networks, local structural properties of the network strongly constrain the
location of the viral epidemic threshold.

The rest of this paper is organized as follows. In Section~\ref{Notation and
Nomenclature}, we review terminology and existing results relating the
dynamical behavior of a virus model with spectral properties of the network
of contacts. In Section \ref{Spectral Analysis from Egonets}, we introduce
an approach, based on algebraic graph theory and convex optimization, to
find upper and lower bounds on the epidemic thresholds from local egonets.
In Subsection~\ref{Sect Moments from Egonets}, we introduce an approach to
related these egonets to the so-called spectral moments of the adjacency
matrix. In Subsection~\ref{Sect Bounds from Moments}, we propose an
optimization framework to derive bounds on the epidemic threshold from a
collection of spectral moments. In Section~\ref{Simulations}, we illustrate
the quality of our approach by computing bounds on the epidemic threshold
for real-world social and communication networks.\medskip

\section{Notation \& Preliminaries\label{Notation and Nomenclature}}

Let $\mathcal{G}=\left( \mathcal{V},\mathcal{E}\right) $ be an undirected,
unweighted graph, where $\mathcal{V}=\left\{ 1,\dots ,n\right\} $ denotes a
set of $n$ nodes and $\mathcal{E}\subseteq \mathcal{V}\times \mathcal{V}$
denotes a set of undirected edges linking them. If $\left( i,j\right) \in 
\mathcal{E}$, we call nodes $i$ and $j$ \emph{adjacent} (or
first-neighbors), which we denote by $i\sim j$. We define the set of
first-neighbors of a node $i$ as $\mathcal{N}_{i}=\{j\in \mathcal{V}:\left(
i,j\right) \in \mathcal{E}\}.$ The \emph{degree} $d_{i}$ of a vertex $i$ is
the number of nodes adjacent to it, i.e., $d_{i}=\left\vert \mathcal{N}%
_{i}\right\vert $. A graph is \emph{weighted} if there is a real number $%
w_{ij}\neq 0$ associated with every edge $\left( i,j\right) \in \mathcal{E}$%
. More formally, a weighted graph $\mathcal{H}$ can be defined as the triad $%
\mathcal{H=}\left( \mathcal{V},\mathcal{E},\mathcal{W}\right) $, where $%
\mathcal{V}$ and $\mathcal{E}$ are the sets of nodes and edges in $\mathcal{H%
}$, and $\mathcal{W=}\left\{ w_{ij}\in \mathbb{R}\backslash \left\{
0\right\} ,\text{ for all }\left( i,j\right) \in \mathcal{E}\right\} $ is
the set of (possibly negative) weights.

The \emph{adjacency matrix} of a simple graph $\mathcal{G}$, denoted by $A_{%
\mathcal{G}}=[a_{ij}]$, is an $n\times n$ symmetric matrix defined
entry-wise as $a_{ij}=1$ if nodes $i$ and $j$ are adjacent, and $a_{ij}=0$
otherwise. For weighted graphs, the entry $a_{ij}$ is equal to the weight $%
w_{ij}$ for $\left( i,j\right) \in \mathcal{E}$; $0$, otherwise. For
undirected graphs, $A_{\mathcal{G}}$ is a symmetric matrix; thus, $A_{%
\mathcal{G}}$ has a full set of $n$ real and orthogonal eigenvectors with
real eigenvalues $\lambda _{1}\geq \lambda _{2}\geq ...\geq \lambda _{n}$.
The largest eigenvalue of $A_{\mathcal{G}}$, $\lambda _{1}$, is called the 
\emph{spectral radius} of $A_{\mathcal{G}}$. If $A$ has nonnegative entries
and is irreducible (i.e., $\mathcal{G}$ is connected), then the
Perron-Frobenius theorem \cite{Mac00} can be used to show that the spectral
radius $\lambda _{1}$ is unique, real, and positive. We also define the $k$%
-th spectral moment of $\mathcal{G}$ as: 
\begin{equation}
m_{k}\left( \mathcal{G}\right) \triangleq \frac{1}{n}\sum_{i=1}^{n}\lambda
_{i}^{k}.  \label{Moments Definition}
\end{equation}

A \emph{walk} of length $k$ from node $i_{1}$ to node $i_{k+1}$ is an
ordered sequence of nodes $\left( i_{1},i_{2},...,i_{k+1}\right) $ such that 
$i_{j}\sim i_{j+1}$ for $j=1,2,...,k$. One says that the walk \emph{touches}
each of the nodes that comprises it. If $i_{1}=i_{k+1}$, then the walk is
closed. A closed walk with no repeated nodes (with the exception of the
first and last nodes) is called a \emph{cycle}. Given a walk $p=\left(
i_{1},i_{2},...,i_{k+1}\right) $ in a weighted graph $\mathcal{H}$, we
define the weight of the walk as, $\omega \left( p\right)
=w_{i_{1}i_{2}}w_{i_{2}i_{3}}...w_{i_{k}i_{k+1}}$.

\subsection{Stochastic Modeling of Viral Spreading}

A wide variety of stochastic models has been proposed in the literature to
study the dynamics of virus spreading processes. In most models, the
steady-state level of infection in the network presents two different
regimes separated by a phase transition taking place at an epidemic
threshold. This epidemic threshold is determined by both the virulence of
the infection and the network topology. A series of papers study the value
of this epidemic threshold as a function of the network structure, in both
random \cite{PV01}--\cite{PJ09} and real topologies \cite{GMT05}--\cite%
{PCFVF11}. A spreading model widely considered in the literature is the
so-called SIS (Susceptible-Infected-Susceptible) model. In this model, each
individual in the network can be in one of two possible states: \emph{%
susceptible} or \emph{infected}. Given an initial set of infected
individuals, the virus propagates through the edges of an undirected graph $%
\mathcal{G}$ at an infection rate $\beta $. Simultaneously, infected nodes
recover at a rate $\delta $, returning back to the susceptible state (see 
\cite{GMT05} for a formal description of this model). In \cite{GMT05}--\cite%
{MOK09}, we find different (and complementary) approaches to find an
expression for the SIS epidemic threshold. In all of these papers, the
authors are able to decouple the effect of the network topology from the
dynamics of individual nodes. On the one hand, the effect of the node
dynamics is completely characterized by the ratio $\tau _{SIS}\triangleq
\delta /\beta $. On the other hand, the effect of the network topology
depends \emph{exclusively} on the largest eigenvalue of the network
adjacency matrix, $\lambda _{1}\left( A_{\mathcal{G}}\right) $, such that if
the threshold condition $\lambda _{1}\left( A\right) <\tau _{SIS}=\delta
/\beta $ is satisfied, a `small' initial infection dies out exponentially
fast \cite{GMT05}--\cite{MOK09}.

Many extensions to the SIS model have been proposed to capture different
characteristics of viral processes, such as permanent or temporal immunity
of a recovered individual, or virus incubation time \cite{DGM08,PCFVF11}. As
shown in \cite{PCFVF11}, the decoupling argument that allows to separate the
role of the network topology from the node dynamics in the SIS model still
holds for a variety of other virus models. Similarly, a `small' initial
infection dies out exponentially fast in these models if the condition $%
\lambda _{1}\left( A_{\mathcal{G}}\right) <\tau _{VM}$ is satisfied, where
the threshold $\tau _{VM}$ measures the virulence of the infection (and is
independent of the network structure). As a bottom line, all of the above
results remark the key role played by the largest eigenvalue of the
adjacency matrix, $\lambda _{1}\left( A_{\mathcal{G}}\right) $, in virus
spreading processes. In particular, the larger $\lambda _{1}\left( A_{%
\mathcal{G}}\right) $, the more efficient a network is to spread a disease
(or a piece of information) virally.

\subsection{Spectral Estimators Based on Random Graphs}

Random network models are currently the primary tool to study the
relationship between local structural properties of a network and its
epidemic threshold. Although many random networks have been proposed in the
literature \cite{NSW01}--\cite{New09}, only random networks including a very
limited amount of structural information are currently amenable to analysis.
The original random graph model is the Erd\"{o}s-R\'{e}nyi graph, denoted by 
$G\left( n,p\right) $, in which each edge in a graph with $n$ nodes is
independently chosen with probability $p$, \cite{ER61}. In this model, the
distribution of degrees in the network follows a Poisson distribution with
expectation $\mathbb{E}[d_{i}]=(n-1)p$. Furthermore, the largest eigenvalue
of its adjacency matrix is almost surely $\lambda _{1}=\left[ 1+o\left(
1\right) \right] np$ (assuming that $np=\Omega \left( \log n\right) $).
Although very interesting from a theoretical point of view, the original
random graph presents very limited modeling capabilities, since the degree
distributions of real-world networks almost never follow a Poisson
distribution.

In order to increase the modeling abilities of random graphs, Chung et al.
proposed in \cite{CLV03} a random graph $G\left( \mathbf{w}\right) $ in
which one can prescribe a desired expected sequence of degrees, $\mathbf{w}%
=\left( w_{1},...,w_{n}\right) $. In this random graph, edges are
independently assigned to each pair of vertices $\left( i,j\right) $ with
probability $\left. w_{i}w_{j}\right/ \sum_{k=1}^{n}w_{k}$. Chung et al.
proved in \cite{CLV03} that if $\left. \sum_{i=1}^{n}w_{i}^{2}\right/
\sum_{j=1}^{n}w_{j}>\sqrt{\max \left\{ w_{i}\right\} }\log n$, then the
largest eigenvalue $\lambda _{1}\left( G\left( \mathbf{w}\right) \right) $
converges almost surely%
\begin{equation}
\lambda _{1}\left( G\left( \mathbf{w}\right) \right) \overset{a.s.}{%
\rightarrow }\left[ 1+o\left( 1\right) \right] \frac{\sum_{i=1}^{n}w_{i}^{2}%
}{\sum_{j=1}^{n}w_{j}},  \label{CL Estimator}
\end{equation}%
for large $n$. Despite its theoretical interest, random graphs with a given
degree distribution are by far not enough to faithfully model the structure
of real complex networks. In particular, it is well-known that the degree
distribution alone is not a sufficient statistic to analyze the performance
of many networks. For example, Alderson et al. introduce in \cite{ALWD05} a
collection of networks, including random graphs, presenting the same degree
distribution and radically different dynamical performance.

Although random graph models with more elaborated structural properties can
be found in the literature \cite{NSW01}--\cite{New09}, these models are
usually hard (if not impossible) to analyze from a spectral point of view.
The source of this intractability is the presence of strong correlations
among the entries of the (random) adjacency matrix. These strong
correlations prevent the resulting random adjacency matrix from being
analytically tractable. In the next section, we present a novel approach to
analyze the effect of local structural properties on the largest eigenvalue
of a network without making use of random graphs.

\section{Spectral Analysis from Egocentric Subnetworks\label{Spectral
Analysis from Egonets}}

In this section, we study the relationship between local structural
properties of a network and its eigenvalue spectrum. In our analysis, we
assume that we do not have access to the complete topology of the network,
due to, for example, privacy and/or security constrains. Instead, we assume
that we are able to access local egocentric views of the network topology.
In this setting, we propose an approach to extract global spectral
information from local structural properties of the network. This spectral
information will be used in Subsection \ref{Sect Bounds from Moments} to
compute upper and lower bounds on the epidemic threshold.

We now provide graph-theoretical and algebraic elements to characterize the
information contained in these egocentric views of the network. Let $\delta
\left( i,j\right) $ denote the \emph{distance} between two nodes $i$ and $j$
(i.e., the minimum length of a walk from $i$ to $j$). By convention, we
assume that $\delta \left( i,i\right) =0$. We define the $r$-th order
neighborhood around node $i$, denoted by $\mathcal{G}_{i,r}=(\mathcal{N}%
_{i,r},\mathcal{E}_{i,r})$,\ as the subgraph $\mathcal{G}_{i,r}\subseteq 
\mathcal{G}$ with node-set $\mathcal{N}_{i,r}\triangleq \left\{ j\in 
\mathcal{V}:\delta \left( i,j\right) \leq r\right\} $, and edge-set $%
\mathcal{E}_{i,r}=\{(v,w)\in \mathcal{E}$ s.t. $v,w\in \mathcal{N}_{i,r}\}$.
Notice that $\mathcal{G}_{i,r}$ provides a graph-theoretical description of
the egocentric view of the network from node $i$ within a radius of $r$
hops. Motivated by this interpretation, we also call $\mathcal{G}_{i,r}$ the 
\emph{egonet} of radius $r$ around node $i$. Egonets can be algebraically
represented via submatrices of the adjacency matrix $A_{\mathcal{G}}$, as
follows. Given a set of $k$ nodes $\mathcal{K}\subseteq \mathcal{V}$, we
denote by $A_{\mathcal{G}}\left( \mathcal{K}\right) $ the $k\times k$
submatrix of $A_{\mathcal{G}}$ formed by selecting the rows and columns of $%
A_{\mathcal{G}}$ indexed by $\mathcal{K}$. In particular, we define the
adjacency submatrix $A_{i,r}\triangleq A_{\mathcal{G}}\left( \mathcal{N}%
_{i,r}\right) $. Notice that $A_{i,r}$ is itself an adjacency matrix
representing the structure of the egonet $\mathcal{G}_{i,r}$. By convention,
we associate the first row and column of the submatrix $A_{i,r}$ with node $%
i\in \mathcal{V}$, which can be done via a simple permutation of the rows
and columns of $A_{i,r}$.\footnote{%
Notice that permuting the rows and columns of the adjacency matrix does not
change the topology of the underlying graph, it simply changes the labels
associated to each node.} For a weighted graph $\mathcal{H}$ with weighted
adjacency matrix $A_{\mathcal{H}}$, we define the weighted egonet $\mathcal{H%
}_{i,r}$ as the weighted graph whose adjacency matrix is $A_{i,r}\triangleq
A_{\mathcal{H}}\left( \mathcal{N}_{i,r}\right) $.

\subsection{Spectral Moments from Local Egonets\label{Sect Moments from
Egonets}}

In this subsection, we derive expressions for the spectral moments of the
adjacency from the knowledge of local egonets using tools from algebraic
graph theory. The following lemma provides an interesting connection between
the number of closed walks in $\mathcal{G}$ (a combinatorial property) and
its spectral moments (an algebraic property) \cite{Big93}:

\medskip

\begin{lemma}
\label{Diagonals as Walks}Let $\mathcal{G}$ be a simple graph with adjacency
matrix $A_{\mathcal{G}}=\left[ a_{ij}\right] $. Then%
\begin{equation*}
\left[ A_{\mathcal{G}}^{k}\right] _{ii}=\left\vert W_{i,k}\right\vert ,
\end{equation*}%
where $W_{i,k}$ is the set of closed walks of length $k$ starting and
finishing at node $i$.
\end{lemma}

\medskip

Using the above result, one can prove the following well-known result in
algebraic graph theory \cite{Big93}:

\medskip

\begin{corollary}
\label{Three Moments}Let $\mathcal{G}$ be a simple graph. Denote by $e$ and $%
\Delta $ the number of edges and triangles in $\mathcal{G}$, respectively.
Then,%
\begin{equation}
m_{1}(A_{\mathcal{G}})=0,~m_{2}(A_{\mathcal{G}})=\frac{2e}{n},~\text{and }%
m_{3}(A_{\mathcal{G}})=\frac{6\Delta }{n}.  \notag
\end{equation}
\end{corollary}

\medskip

We can generalize Lemma \ref{Diagonals as Walks} to weighted graphs as
follows:

\medskip

\begin{lemma}
\label{Moments from Walks}Let $\mathcal{H=}\left( \mathcal{V},\mathcal{E},%
\mathcal{W}\right) $ be a weighted graph with weighted adjacency matrix $A_{%
\mathcal{H}}$. Then, 
\begin{equation*}
\left[ A_{\mathcal{H}}^{k}\right] _{ii}=\sum_{p\in P_{k,i}}\omega \left(
p\right) ,
\end{equation*}%
where $P_{k,i}$ is the set of closed walks of length $k$ from $v_{i}$ to
itself in $\mathcal{H}$.
\end{lemma}

\medskip

\begin{proof}
By recursively applying the multiplication rule for matrices, we have the
following expansion%
\begin{equation}
\left[ A_{\mathcal{H}}^{k}\right] _{ii}=\sum_{i=1}^{n}\sum_{i_{2}=1}^{n}%
\cdots \sum_{i_{k}=1}^{n}w_{i,i_{2}}w_{i_{2,}i_{3}}\cdots ~w_{i_{k},i}.
\label{Multiple Summations}
\end{equation}%
Using the graph-theoretic nomenclature introduced in Section \ref{Notation
and Nomenclature}, we have that $w_{i,i_{2}}w_{i_{2,}i_{3}}...w_{i_{k},i}=%
\omega \left( p\right) $, for $p=\left(
v_{i},v_{i_{2}},v_{i_{3}},...,v_{i_{k}},v_{i}\right) $. Hence, the
summations in (\ref{Multiple Summations}) can be written as $\left[ W_{%
\mathcal{H}}^{k}\right] _{ii}=\sum_{1\leq i,i_{2},...,i_{k}\leq n}\omega
\left( p\right) $. Finally, the set of closed walks $p=\left(
v_{i},v_{i_{2}},v_{i_{3}},...,v_{i_{k}},v_{i}\right) $\ with indices $1\leq
i,i_{2},...,i_{k}\leq n$ is equal to the set of closed walks of length $k$
from $v_{i}$ to itself in $\mathcal{H}$ (which we have denoted by $P_{k,i}$
in the statement of the Proposition).
\end{proof}

\begin{figure}[tbp]
\centering\includegraphics[width=0.65\linewidth]{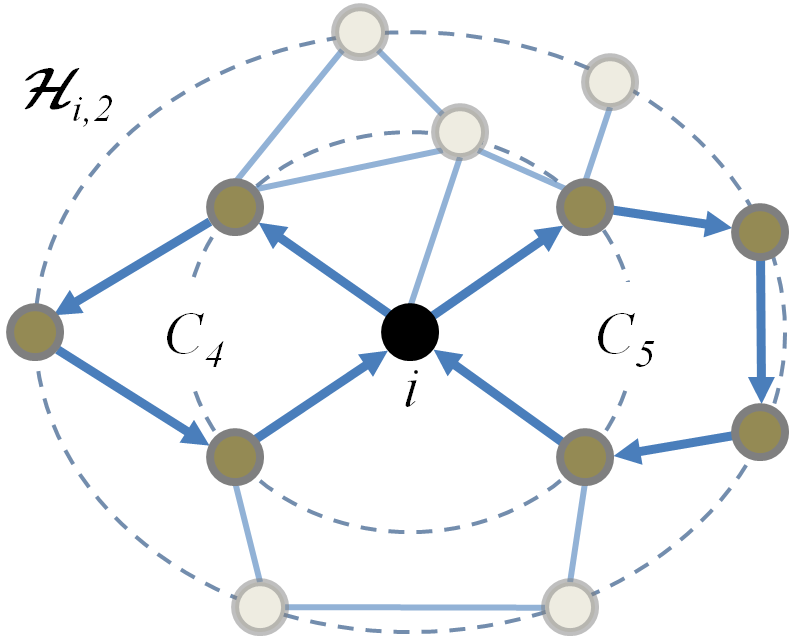}
\caption{Cycles $\mathcal{C}_{4}$ and $\mathcal{C}_{5}$, of lengths $4$ and $%
5\,$, in a neighborhood of radius $2$ around node \thinspace $i\,$.}
\end{figure}

\medskip

Using Lemma \ref{Moments from Walks}, we can extend Lemma \ref{Three Moments}
to higher-order moments of weighted graphs as follows:

\medskip

\begin{theorem}
\label{Moments from Subgraphs}Consider a weighted, undirected graph $%
\mathcal{H}$ with adjacency matrix $A_{\mathcal{H}}$. Let $A_{i,r}$ be the
(weighted) adjacency matrix of the egonet of radius $r$ around node $i$.
Then, for a given $r$, the spectral moments of $A_{\mathcal{H}}$ can be
written as%
\begin{equation}
m_{k}\left( A_{\mathcal{H}}\right) =\frac{1}{n}\sum_{i=1}^{n}\left[
A_{i,r}^{k}\right] _{11},  \label{Moments as Sum of Traces}
\end{equation}%
for $k\leq 2r+1$.
\end{theorem}

\medskip

\begin{proof}
Since the trace of a matrix is the sum of its eigenvalues, we can expand the 
$k$-th spectral moment of the adjacency matrix as follows:%
\begin{equation}
m_{k}\left( A_{\mathcal{H}}\right) =\frac{1}{n}\text{Trace}\left( A_{%
\mathcal{H}}^{k}\right) =\frac{1}{n}\sum_{i=1}^{n}\left[ A_{\mathcal{H}}^{k}%
\right] _{ii}.  \label{Moments as Walks in Laplacian Graph}
\end{equation}%
From Lemma \ref{Moments from Walks}, we have that $\left[ A_{\mathcal{H}}^{k}%
\right] _{ii}=\sum_{p\in P_{k,i}}\omega \left( p\right) $. Notice that for a
fixed value of $k$, closed walks of length $k$ in $\mathcal{H}$ starting at
node $i$ can only touch nodes within a certain distance $r\left( k\right) $
of $i$, where $r\left( k\right) $ is a function of $k$. In particular, for $%
k $ even (resp. odd), a closed walk of length $k$ starting at node $i$ can
only touch nodes at most $k/2$ $\,$(resp. $\left\lfloor k/2\right\rfloor $)
hops away from $i$ (see Fig. 1). Therefore, closed walks of length $k$
starting at $i$ are always contained within the neighborhood of radius $%
\left\lfloor k/2\right\rfloor $. In other words, the egonet $\mathcal{H}%
_{i,r}$ of radius $r$ contains all closed walks of length up to $2r+1$
starting at node $i$. We can count these walks by applying Lemma \ref%
{Moments from Walks} to the local adjacency matrix $A_{i,r}$. In particular, 
$\sum_{p\in P_{k,i}}\omega \left( p\right) $ is equal to $\left[ A_{i,r}^{k}%
\right] _{11}$ (since, by convention, node $1$ in the local egonet $\mathcal{%
H}_{i,r}$ corresponds to node $i$ in the graph $\mathcal{H}$). Therefore,
for $k\leq 2r+1$, we have that 
\begin{equation}
\left[ A_{i,r}^{k}\right] _{11}=\sum_{p\in P_{k,i}}\omega \left( p\right) =%
\left[ A_{\mathcal{H}}^{k}\right] _{ii}.
\label{Walks as Diagonal of Sublaplacian}
\end{equation}%
Then, substituting (\ref{Walks as Diagonal of Sublaplacian}) into (\ref%
{Moments as Walks in Laplacian Graph}), we obtain the statement of our
Theorem.
\end{proof}

\medskip

\begin{remark}
The above theorem allows us to compute a truncated sequence of spectral
moments $\left\{ m_{k}\left( A_{\mathcal{H}}\right) \text{, }k\leq
2r+1\right\} $, given a collection of local egonets of radius $r$, $\left\{ 
\mathcal{H}_{i,r}\text{, }i\in \mathcal{V}\right\} $. According to (\ref%
{Moments as Sum of Traces}), we can compute the $k$-th spectral moment by
simply averaging the quantities $\left[ A_{i,r}^{k}\right] _{11}$, $%
i=1,...,n $. For a fixed $k$, each value $\left[ A_{i,r}^{k}\right] _{11}$, $%
i=1,\dots,n$, can be computed in time $O\left( \left\vert \mathcal{N}%
_{i,r}\right\vert ^{3}\right) $, where $\left\vert \mathcal{N}%
_{i,r}\right\vert $ is the number of nodes in the local egonet $\mathcal{H}%
_{i,r}$. The sparse structure of most real networks implies that $\left\vert 
\mathcal{N}_{i,r}\right\vert \ll n$ (for moderate values of $r$). In
particular, if $\left\vert \mathcal{N}_{i,r}\right\vert =o\left(
n^{\epsilon}\right) $ for any $\epsilon>0$, we can compute the $k$-th
spectral moments in quasi-linear time (with respect to the size of the
network) using (\ref{Moments as Sum of Traces}). This result provides a
clear computational advantage compared to computing the spectral moments via
an explicit eigenvalue decomposition, which can be prohibitively expensive
to compute for large complex networks.
\end{remark}


\subsection{SDP-Based Bounds on the Spectral Radius\label{Sect Bounds from
Moments}}

Using Theorem \ref{Moments from Subgraphs}, we can compute a truncated
sequence of the spectral moments of a network $\mathcal{H}$, $(m_{1}\left(
A_{\mathcal{H}}\right) ,m_{2}\left( A_{\mathcal{H}}\right)
,...,m_{2r+1}\left( A_{\mathcal{H}}\right) )$, from a set of local egonets
of radius $r$, $\left\{ \mathcal{H}_{i,r}\text{, }i\in \mathcal{V}\right\} $%
. We now present a convex optimization framework to extract information
about the largest eigenvalue of the adjacency matrix, $\lambda _{1}\left( A_{%
\mathcal{H}}\right) $, from this sequence of moments. We can state the
problem solved in this subsection as follows:

\medskip

\begin{problem}
\label{Bound from Moments problem}Given a truncated sequence of spectral
moments of a weighted, undirected graph $\mathcal{H}$, $\mathbf{m}%
_{2r+1}=\left( m_{0},m_{1},...,m_{2r+1}\right) $, find tight upper and lower
bounds on the largest eigenvalue $\lambda _{1}\left( A_{\mathcal{H}}\right) $%
.
\end{problem}

\medskip

Our approach is based on a probabilistic interpretation of the eigenvalue
spectrum of a given network. To present our approach, we first need to
introduce some concepts:

\medskip

\begin{definition}
\label{Spectral density}Given a weighted, undirected graph $\mathcal{H}$
with (real) eigenvalues $\lambda _{1},...,\lambda _{n}$, the \emph{spectral
density} of $\mathcal{H}$ is defined as,%
\begin{equation}
\mu _{\mathcal{H}}\left( x\right) \triangleq \frac{1}{n}\sum_{i=1}^{n}\delta
\left( x-\lambda _{i}\right) ,  \label{Spectral Measure}
\end{equation}%
where $\delta \left( \cdot \right) $ is the Dirac delta function.
\end{definition}

\medskip

The spectral density can be interpreted as a discrete probability density
function with \emph{support\footnote{%
Recall that the support of a finite Borel measure $\mu $ on $\mathbb{R}$,
denoted by $supp\left( \mu \right) $, is the smallest closed set $B$ such
that $\mu \left( \mathbb{R}\backslash B\right) =0$.}} on the set of
eigenvalues $\left\{ \lambda _{i},\text{ }i=1...n\right\} $. Let us consider
a discrete random variable $X$ whose probability density function is $\mu _{%
\mathcal{H}}$. The moments of this random variable satisfy the following:

\medskip

\begin{lemma}
\label{Spectral moments}The moments of a r.v. $X\sim \mu _{\mathcal{H}}$ are
equal to the spectral moments of $A_{\mathcal{H}}$, i.e.,%
\begin{equation*}
\mathbb{E}_{\mu _{\mathcal{H}}}\left( X^{k}\right) =m_{k}\left( A_{\mathcal{H%
}}\right) ,
\end{equation*}%
for all $k\geq 0$.
\end{lemma}

\medskip

\begin{proof}
For all $k\geq 0$, we have the following:%
\begin{align*}
\mathbb{E}_{\mu _{\mathcal{H}}}\left( X^{k}\right) & =\int_{\mathbb{R}%
}x^{k}\mu _{\mathcal{H}}\left( x\right) dx \\
& =\frac{1}{n}\sum_{i=1}^{n}\int_{\mathbb{R}}x^{k}\delta \left( x-\lambda
_{i}\right) dx \\
& =\frac{1}{n}\sum_{i=1}^{n}\lambda _{i}^{k}=m_{k}\left( A_{\mathcal{G}%
}\right) .
\end{align*}
\end{proof}

\medskip

We now present a convex optimization framework that allows us to find bounds
on the endpoints of the smallest interval $\left[ a,b\right] $ containing
the support of a generic random variable $X\sim \mu $ given a sequence of
moments $\left( M_{0},M_{1},...,M_{2r+1}\right) $, where $M_{k}\triangleq
\int x^{k}d\mu $. Subsequently, we shall apply these results to find bounds
on $\lambda _{1}\left( A_{\mathcal{H}}\right) $. Our formulation is based on
the following matrices:

\medskip

\begin{definition}
\label{Hankel Moment Matrices}Given a sequence of moments $\mathbf{M}%
_{2r+1}=\left( M_{0},M_{1},...,M_{2r+1}\right) $, let $H_{2r}\left( \mathbf{M%
}_{2r+1}\right) $ and $H_{2r+1}\left( \mathbf{M}_{2r+1}\right) \in \mathbb{R}%
^{\left( r+1\right) \times \left( r+1\right) }$ be the Hankel matrices
defined by\footnote{%
For simplicity in the notation, we shall omit the argument $\mathbf{M}%
_{2r+1} $ whenever clear from the context.}:%
\begin{equation}
\left[ H_{2r}\right] _{ij}\triangleq M_{i+j-2}\text{ and }\left[ H_{2r+1}%
\right] _{ij}\triangleq M_{i+j-1}.  \label{Moment Matrices}
\end{equation}%
The above matrices are called the \emph{moment matrices} associated with the
sequence $\mathbf{M}_{2r+1}$.
\end{definition}

\medskip

In general, an arbitrary sequence of numbers $\left(
N_{0},N_{1},...,N_{k}\right) $ may not have a representing measure $\mu $
such that $\int x^{r}d\mu =N_{r}$, for $0\leq r\leq k$. A sequence of
numbers $\mathbf{N}_{k}=\left( N_{0},N_{1},...,N_{k}\right) $ is said to be 
\emph{feasible} in $\Omega \subseteq \mathbb{R}$ if there exists a measure $%
\mu $ with support contained in $\Omega $ whose moments match those in the
sequence $\mathbf{N}_{k}$. The problem of deciding whether or not a sequence
of numbers is feasible in $\Omega $ is called the \emph{classical moment
problem} in analysis \cite{ST43}. For univariate distributions, necessary
and sufficient conditions for feasibility can be given in terms of certain
Hankel matrices being positive semidefinite\footnote{%
The notation $A\succeq 0$ means that the matrix $A$ is positive semidefinite.%
}, as follows \cite{Las09}:

\medskip

\begin{theorem}
{\cite[Theorem 3.2]{Las09}} \label{Feasibility Results}Let $\mathbf{M}%
_{2r+1}=\left( M_{0},M_{1},...,M_{2r+1}\right) \in \mathbb{R}^{2r+2}$. Then,

\begin{description}
\item[(a)] The sequence $\mathbf{M}_{2r+1}$ corresponds to a sequence of
moments feasible in $\Omega =\mathbb{R}$ if and only if $H_{2r}\succeq 0$.

\item[(b)] The sequence $\mathbf{M}_{2r+1}$ is feasible in $\Omega =\left[
a,\infty \right) $ if and only if 
\begin{equation*}
H_{2r}\succeq 0\text{ and }H_{2r+1}-aH_{2r}\succeq 0.
\end{equation*}

\item[(c)] The sequence $\mathbf{M}_{2r+1}$ is feasible in $\Omega =\left(
-\infty ,b\right] $ if and only if%
\begin{equation*}
H_{2r}\succeq 0\text{ and }bH_{2r}-H_{2r+1}\succeq 0.
\end{equation*}
\end{description}
\end{theorem}

\medskip


\medskip

Using Theorem \ref{Feasibility Results}, we have the following result:

\medskip

\begin{theorem}
\label{Main Theorem for general densities}Let $\mu $ be a probability
density function on $\mathbb{R}$ with associated sequence of moments $%
\mathbf{M}_{2r+1}=\left( M_{0},M_{1},...,M_{2r+1}\right) $, all finite, and
let $\left[ a,b\right] $ be the smallest interval which contains the support
of $\mu $. Then, $b\geq \beta ^{\ast }\left( \mathbf{M}_{2r+1}\right) $,
where%
\begin{equation}
\begin{array}{rrl}
\beta _{r}^{\ast }\left( \mathbf{M}_{2r+1}\right) \triangleq & \min_{x} & x
\\ 
& \text{s.t.} & H_{2r}\succeq 0, \\ 
&  & x~H_{2r}-H_{2r+1}\succeq 0.%
\end{array}
\label{SDP generic bound}
\end{equation}
\end{theorem}

\medskip

\begin{proof}
Since $\mathbf{M}_{2r+1}$ is the moment sequence of a probability density
function $\mu $ with support on $\left[ a,b\right] \subset \left( -\infty ,b%
\right] $, we have from Theorem \ref{Feasibility Results}.(\emph{c}) that $%
\mathbf{M}_{2r+1}$ satisfy $H_{2r}\succeq 0$ and $bH_{2r}-H_{2r+1}\succeq 0$%
. Since $\beta ^{\ast }\left( \mathbf{M}_{2r+1}\right) $ is, by definition,
the minimum value of $x$ such that $H_{2r}\succeq 0$ and $%
xH_{2r}-H_{2r+1}\succeq 0$, we have that $\beta ^{\ast }\left( \mathbf{M}%
_{2r+1}\right) \leq b$.
\end{proof}

\medskip

\begin{remark}
Observe that, for a given sequence of moments $\mathbf{M}_{2r+1}$, the
entries of $xH_{2r}-H_{2r+1}$ depend affinely on the variable $x$. Then $%
\beta ^{\ast }\left( \mathbf{m}_{2r+1}\right) $ is the solutions to a
semidefinite program\footnote{%
A semidefinite program is a convex optimization problem that can be solved
in time polynomial in the input size of the problem; see e.g. \cite{VB96}.}
(SDP) in one variable. Hence, $\beta ^{\ast }\left( \mathbf{M}_{2r+1}\right) 
$ can be efficiently computed using standard optimization software, e.g. 
\cite{CVX}, from a truncated sequence of moments.
\end{remark}

Applying Theorem \ref{Main Theorem for general densities} to the spectral
density $\mu _{\mathcal{H}}$ of a given graph $\mathcal{H}$ with spectral
moments $\left( m_{0},m_{1},...,m_{2r+1}\right) $, we can find a lower bound
on its largest eigenvalue, $\lambda _{1}\left( A_{\mathcal{H}}\right) $, as
follows:

\medskip

\begin{theorem}
\label{Main Theorem for Spectral Densities}Let $\mathcal{H}$ be a weighted,
undirected graph with (real) eigenvalues $\lambda _{1}\geq ...\geq \lambda
_{n}$. Then, given a truncated sequence of the spectral moments of $\mathcal{%
H}$, $\mathbf{m}_{2r+1}=\left( m_{0},m_{1},...,m_{2r+1}\right) $, we have
that%
\begin{equation}
\lambda _{1}\left( A_{\mathcal{H}}\right) \geq \beta _{r}^{\ast }\left( 
\mathbf{m}_{2r+1}\right) ,  \label{SDP spectral bound}
\end{equation}%
(where $\beta _{r}^{\ast }\left( \mathbf{m}_{2r+1}\right) $ is the solution
to the SDP in (\ref{SDP generic bound})).
\end{theorem}

\medskip

\begin{proof}
Let us consider the spectral density of $\mathcal{H}$, $\mu _{\mathcal{H}}$,
in Definition \ref{Spectral density}. According to Lemma \ref{Spectral
moments}, the density $\mu _{\mathcal{H}}$ has associated moments $\mathbf{m}%
_{2r+1}$. Also, the smallest interval which contains the support of $\mu _{%
\mathcal{H}}$ is $\left[ a,b\right] =\left[ \lambda _{n},\lambda _{1}\right] 
$. Therefore, applying Theorem \ref{Main Theorem for general densities} to $%
\mu _{\mathcal{H}}$, we obtain that $\beta _{r}^{\ast }\left( \mathbf{m}%
_{2r+1}\right) \leq b=\lambda _{1}$.
\end{proof}

\medskip

Furthermore, for $r=1$, we can analytically solve the SDP in (\ref{SDP
generic bound}) to derive a closed-form solution for $\beta _{1}^{\ast
}\left( \mathbf{m}_{3}\right) $, as follows:

\medskip

\begin{corollary}
\label{Three Moments Bound}Let $\mathcal{G}$ be a simple graph with
adjacency matrix $A_{\mathcal{G}}$. Denote by $n$, $e$, and $\Delta $ the
number of nodes, edges, and triangles in $\mathcal{G}$, respectively. Then,%
\begin{equation}
\lambda _{1}\left( A_{\mathcal{G}}\right) \geq \frac{3\Delta +\sqrt{9\Delta
^{2}+8e^{3}/n}}{2e}.  \label{Bound from three}
\end{equation}
\end{corollary}

\medskip

\begin{proof}
In the Appendix.
\end{proof}

\medskip

Using the optimization framework presented above, we can also compute upper
bounds on the spectral radius of $\mathcal{H}$ from a sequence of its
spectral moments, as follows. In this case, our formulation is based on the
following set of Hankel matrices:

\medskip

\begin{definition}
Given a weighted, undirected graph $\mathcal{H}$ with $n$ nodes and spectral
moments $\mathbf{m}_{2r+1}=\left( m_{0},m_{1},...,m_{2r+1}\right) $, let $%
T_{2r}\left( y;\mathbf{m}_{2r+1},n\right) $ and $T_{2r+1}\left( y;\mathbf{m}%
_{2r+1},n\right) \in \mathbb{R}^{\left( r+1\right) \times \left( r+1\right)
} $ be the Hankel matrices defined by\footnote{%
We shall omit the arguments from $T_{2r}$ and $T_{2r+1}$ whenever clear from
the context.}:%
\begin{eqnarray}
\left[ T_{2r}\right] _{ij} &\triangleq &\frac{n}{n-1}m_{i+j-2}-\frac{1}{n-1}%
y^{i+j-2},  \label{Bulk Hankel Matrices} \\
\left[ T_{2r+1}\right] _{ij} &\triangleq &\frac{n}{n-1}m_{i+j-1}-\frac{1}{n-1%
}y^{i+j-1}.  \notag
\end{eqnarray}
\end{definition}

\medskip

Given a sequence of spectral moments, we can compute upper bounds on the
largest eigenvalue $\lambda _{1}\left( A_{\mathcal{H}}\right) $ using the
following result:

\medskip

\begin{theorem}
\label{Upper Bound for Spectral Densities}Let $\mathcal{H}$ be a weighted,
undirected graph with (real) eigenvalues $\lambda _{1}\geq ...\geq \lambda
_{n}$. Then, given a truncated sequence of its spectral moments $\mathbf{m}%
_{2r+1}=\left( m_{0},m_{1},...,m_{2r+1}\right) $, we have that 
\begin{equation*}
\lambda _{1}\leq \delta _{r}^{\ast }\left( \mathbf{m}_{2r+1},n\right) ,
\end{equation*}%
where 
\begin{equation}
\begin{array}{rrl}
\delta _{r}^{\ast }\left( \mathbf{m}_{2r+1},n\right) \triangleq & \max_{y} & 
y \\ 
& \text{s.t.} & T_{2r}\succeq 0, \\ 
&  & yT_{2r}-T_{2r+1}\succeq 0, \\ 
&  & T_{2r+1}+yT_{2r}\succeq 0.%
\end{array}
\label{Upper Bound from Moments}
\end{equation}
\end{theorem}

\medskip

\begin{proof}
Let us define the \emph{bulk} of the spectrum as the set of eigenvalues $%
\left\{ \lambda _{2},...,\lambda _{n}\right\} $, and the \emph{bulk spectral
density} as the probability density function:%
\begin{equation*}
\tilde{\mu}_{\mathcal{H}}\triangleq \frac{1}{n-1}\sum_{i=2}^{n}\delta \left(
x-\lambda _{i}\right) .
\end{equation*}%
We also define the \emph{bulk spectral moments} as the moments of the bulk
spectral density, which satisfy:%
\begin{align*}
\tilde{m}_{k}\left( A_{\mathcal{H}}\right) & \triangleq \int_{\mathbb{R}%
}x^{k}\tilde{\mu}_{\mathcal{H}}\left( x\right) dx \\
& =\frac{1}{n-1}\sum_{i=2}^{n}\int_{\mathbb{R}}x^{k}\delta \left( x-\lambda
_{i}\right) dx \\
& =\frac{1}{n-1}\sum_{i=1}^{n}\lambda _{i}^{k}-\frac{1}{n-1}\lambda _{1}^k \\
& =\frac{n}{n-1}m_{k}\left( A_{\mathcal{H}}\right) -\frac{1}{n-1}\lambda
_{1}^k.
\end{align*}%
Therefore, the moment matrices associated to the sequence of bulk spectral
moments $\mathbf{\tilde{m}}_{2r+1}=\left( \tilde{m}_{0},\tilde{m}_{1},...,%
\tilde{m}_{2r+1}\right) $, satisfy%
\begin{equation}
H_{s}\left( \mathbf{\tilde{m}}_{2r+1}\right) =T_{s}\left( \lambda _{1};%
\mathbf{m}_{2r+1},n\right) ,  \label{H as T}
\end{equation}%
for $s\in \left\{ 2r,2r+1\right\} $, where $H_{s}$ and $T_{s}$ were defined
in (\ref{Moment Matrices}) and (\ref{Bulk Hankel Matrices}), respectively.

Since $\left\vert \lambda _{i}\right\vert \leq \lambda _{1}$ for $i\geq 2$,
the support of the bulk spectral density $\tilde{\mu}_{\mathcal{H}}$ is
contained in the interval $\left[ -\lambda _{1},\lambda _{1}\right] $.
Hence, according to Theorems \ref{Feasibility Results}.(\emph{b})-(\emph{c}%
), the sequence of bulk spectral moments $\mathbf{\tilde{m}}_{2r+1}$ must
satisfy:%
\begin{equation*}
\begin{array}{l}
T_{2r}\left( \lambda _{1};\mathbf{m}_{2r+1},n\right) \succeq 0, \\ 
\lambda _{1}T_{2r}\left( \lambda _{1};\mathbf{m}_{2r+1},n\right)
-T_{2r+1}\left( \lambda _{1};\mathbf{m}_{2r+1},n\right) \succeq 0, \\ 
T_{2r+1}\left( \lambda _{1};\mathbf{m}_{2r+1},n\right) +\lambda
_{1}T_{2r}\left( \lambda _{1};\mathbf{m}_{2r+1},n\right) \succeq 0.%
\end{array}%
\end{equation*}%
Since $\delta _{r}^{\ast }\left( \mathbf{m}_{2r+1},n\right) $ is, by
definition, the maximum value of $y$ satisfying the constrains in (\ref%
{Upper Bound from Moments}), we have that $\delta _{r}^{\ast }\left( \mathbf{%
m}_{2r+1},n\right) \geq \lambda _{1}\left( A_{\mathcal{H}}\right) $.
\end{proof}

\medskip

\begin{remark}
The optimization program in (\ref{Upper Bound from Moments}) is not an SDP,
since the entries of the matrices $T_{2r}\left( y;\mathbf{m}_{2r+1},n\right) 
$ and $T_{2r+1}\left( y;\mathbf{m}_{2r+1},n\right) $\ are not affine
functions, but higher-order polynomials, in $y$. Nevertheless, the program
can be cast into a convex optimization program, as follows. For the matrices
in (\ref{Upper Bound from Moments}) to be positive semidefinite, all their
principal minors must be nonnegative, where each minor is a polynomial in $y$%
. In other words, positive semidefiniteness of the matrices in (\ref{Upper
Bound from Moments}) is equivalent to a collection of polynomials in $y$
being nonnegative. Hence, we can substitute the semidefinite constrains in (%
\ref{Upper Bound from Moments}) by a collection of polynomials in $y$ being
nonnegative. The resulting optimization problem is a Sum-Of-Squares (SOS)
program \cite{Par00}, which is a type of convex program that can be
efficiently solved using off-the-shelf software \cite{PPSP04}.
\end{remark}

\medskip

In summary, using Theorems \ref{Moments from Subgraphs}, \ref{Main Theorem
for Spectral Densities}, and \ref{Upper Bound for Spectral Densities}, we
can compute upper and lower bounds on the largest eigenvalue of a weighted,
undirected network, $\lambda _{1}\left( A_{\mathcal{H}}\right) $, from the
set of local egonets with radius $r$, as follows: (\emph{1}) Using (\ref%
{Moments as Sum of Traces}), compute the truncated sequence of moments $%
\left( m_{0},m_{1},...,m_{2r+1}\right) $ from the set of egonets, $\{A_{i,r}$%
, $i\in \mathcal{V\}}$, and (\emph{2}) using Theorems \ref{Main Theorem for
Spectral Densities} and \ref{Upper Bound for Spectral Densities}, compute
the upper and lower bounds, $\delta _{r}^{\ast }\left( \mathbf{m}%
_{2r+1},n\right) $ and $\beta _{r}^{\ast }\left( \mathbf{m}_{2r+1}\right) $,
respectively.

\medskip

\section{Numerical Simulations\label{Simulations}}

In this section, we analyze real data from several social and communication
networks to numerically verify the tightness of our bounds. In our first set
of simulations, we study a regional network of Facebook that spans $63,731$
users (nodes) connected by $817,090$ friendships (edges) \cite{VMCG09}. In
order to corroborate our results in different network topologies, we extract
multiple medium-size social subgraphs by running a Breath-First Search (BFS)
around a collection of starting nodes in the Facebook graph. Each BFS
induces a social subgraph spanning all nodes 2 hops away from a starting
node. As a result, we generate a set of 100 different social subgraphs, $%
\mathbf{G}=\{G_{i}\}_{i\leq 100}$, centered around 100 randomly chosen
nodes. For each social subgraph $G_{i}\in \mathbf{G}$, we compute its first
five spectral moments $\mathbf{m}_{5}\left( G_{i}\right) =\left( m_{1}\left(
G_{i}\right) ,...,m_{5}\left( G_{i}\right) \right) $ and use Theorems \ref%
{Main Theorem for Spectral Densities} and \ref{Upper Bound for Spectral
Densities} to compute lower and upper bounds on the spectral radius, $\beta
_{2}^{\ast }\left( G_{i}\right) =\beta _{2}^{\ast }\left( \mathbf{m}%
_{5}\left( G_{i}\right) \right) $ and $\delta _{2}^{\ast }\left(
G_{i}\right) =\delta _{2}^{\ast }\left( \mathbf{m}_{5}\left( G_{i}\right)
,n_{i}\right) $, where $n_{i}$ is the size of $G_{i}$. Since we have access
to the complete network topology, we can also numerically compute the exact
value of the largest eigenvalue $\lambda _{1}\left( G_{i}\right) $, for
comparison purposes. It is worth remarking that, in many real applications,
we do not have access to the complete network topology, due to privacy
and/or security constrains; therefore, we would not be able to compute the
exact value of $\lambda _{1}$. It is in those cases when our approach is
most useful.

Fig. 2 represents a scatter plot where each red circle above the dashed
diagonal line has coordinates $\left( \lambda _{1}\left( G_{i}\right)
,\delta _{2}^{\ast }\left( G_{i}\right) \right) $, and each blue circle
below the dashed diagonal line has coordinates $\left( \lambda _{1}\left(
G_{i}\right) ,\beta _{2}^{\ast }\left( G_{i}\right) \right) $, for all $%
G_{i}\in \mathbf{G}$. We have also included a black line connecting every
pair of circles associated to the same subgraph $G_{i}$. This black line
represents the interval of possible values in which the largest eigenvalue, $%
\lambda _{1}\left( G_{i}\right) $, must lie. (Notice how the dashed diagonal
line cut through all those segments.) For all the social subnetworks in $%
\mathbf{G}$, the spectral radii $\lambda _{1}\left( G_{i}\right) $ are
remarkably close to the theoretical bounds $\beta _{2}^{\ast }\left(
G_{i}\right) $ and $\delta _{2}^{\ast }\left( G_{i}\right) $. In other
words, in our collection of social subgraphs, local structural properties of
the network strongly constrain the location of the largest eigenvalue, and
consequently the ability of a social network to disseminate information
virally.

Our bounds are also tight for other important social and communication
networks. In the following, we the compare the values of $\beta _{2}^{\ast }$
and $\delta _{2}^{\ast }$ with the largest eigenvalue $\lambda _{1}$ of an
e-mail and an Internet network:

\medskip

\begin{example}[Enron e-mail network]
In this example we consider a subgraph of the Enron e-mail communication
network \cite{BY04}. Nodes of the network are e-mail addresses and the
network contains an edge $\left( i,j\right) $ if $i$ sent at least one
e-mail to $j$ (or vice versa). The total size of the network is $36,692$
nodes, which is too large for us to manage computationally. In order to
compare our bounds with the exact value of the largest eigenvalue, we
analyze a subgraph obtained by a BFS of depth 2 around a randomly chosen
node. The resulting subgraph has $n=3,215$ nodes and $e=36,537$ edges. We
also compute the value of its largest eigenvalue to be $\lambda _{1}=95.18$.
Using (\ref{Moments as Sum of Traces}), we have the following values for the
first five spectral moments of the adjacency matrix: $m_{1}=0$, $m_{2}=22.47$%
, $m_{3}=394.7$, $m_{4}=33,491$, and $m_{5}=2,603,200$. From (\ref{SDP
generic bound}) and (\ref{Upper Bound from Moments}), we obtain the
following upper and lower bounds on the largest eigenvalue: $\beta
_{2}^{\ast }=78.53<\lambda _{1}<98.74=\delta _{2}^{\ast }$. Notice that the
numerical value of $\lambda _{1}$ is remarkably close to the upper bound $%
\delta _{2}^{\ast }$. Since the spectral radius measures the ability of a
network to spread information virally, our numerical results indicate that
the e-mail network spreads information very efficiently given the structural
constrains imposed by the local egonets. We can also compare our bounds with
the estimator in (\ref{CL Estimator}), corresponding to a random network
with the same degree distribution. The value of the estimator is equal to $%
\tilde{\lambda}_{1}=124.57$, which is looser than our bounds.
\end{example}

\medskip

\begin{example}[AS-Skitter Internet network]
In this example, we consider a subgraph of the Internet network at the
Autonomous Systems (AS) level. The network topology was obtained from the
Skitter data collection in CAIDA \cite{CAIDA05}. Our subgraph was obtained
from the complete AS graph using a BFS of depth 2 around a random node. The
resulting subgraph has $n=2,248$ nodes, $e=20,648$ edges, and its largest
eigenvalue at $\lambda _{1}=91.3$. The spectral moments of its adjacency
matrix are $m_{1}=0$, $m_{2}=18.37$, $m_{3}=341.1$, $m_{4}=40,001$, and $%
m_{5}=2,777,018$. The resulting bounds from (\ref{SDP generic bound}) and (%
\ref{Upper Bound from Moments}) are $\beta _{2}^{\ast }=74.72<\lambda
_{1}<93.94=\delta _{2}^{\ast }$. Notice how, the largest eigenvalue is again
remarkably close to the upper bound, indicating that the network is able to
spread information efficiently, given its local structural constrains. In
this case, the estimator based on random networks produces a value of $%
\tilde{\lambda}_{1}=219.1$, which is very loose. Therefore, using random
networks to analyze spreading processes in the Internet graph can be
misleading.
\end{example}

\begin{figure}[t]
\centering
\includegraphics[width=0.43\textwidth]{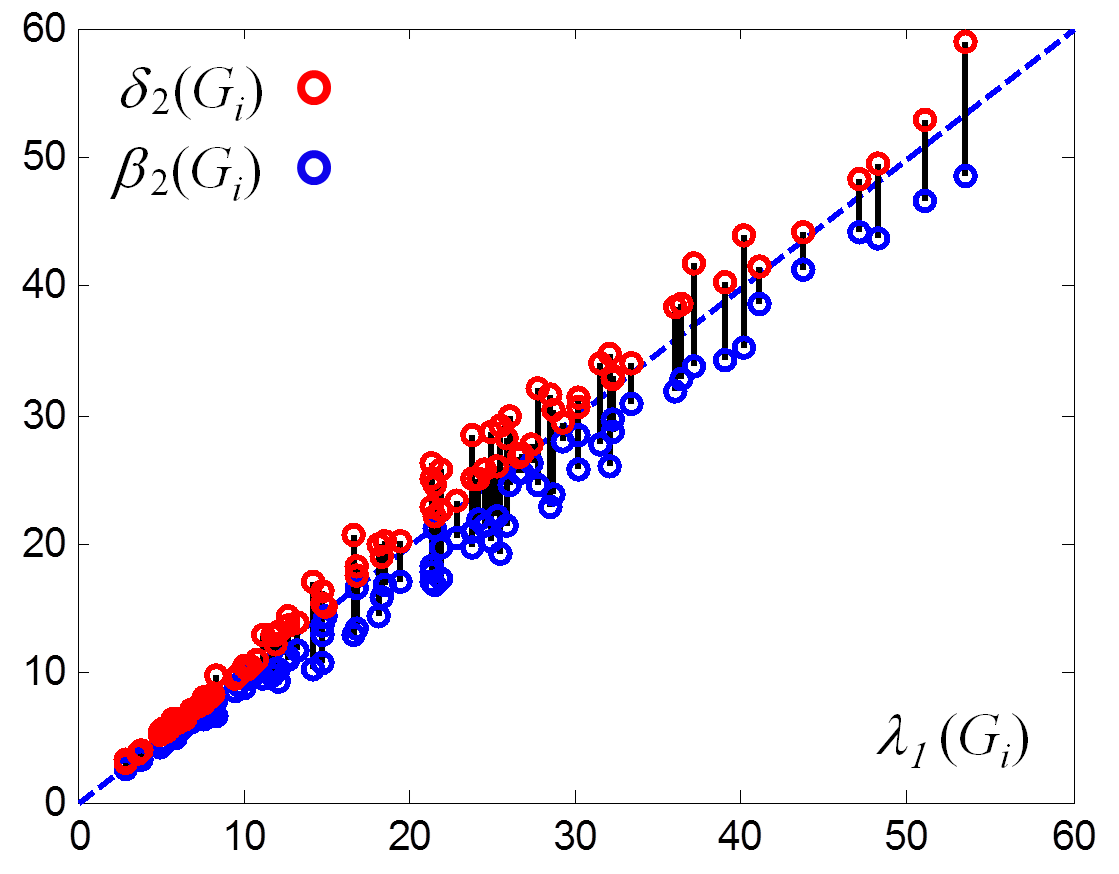}
\caption{Scatter plot of the spectral radius, $\protect\lambda _{1}\left(
G_{i}\right) $, versus the lower bound $\protect\beta _{2}\left(
G_{i}\right) $ (blue circles) and the upper bound $\protect\delta _{2}\left(
G_{i}\right) $ (red circles), where each point is associated with one of the 
$100$ social subgraphs considered in our experiments.}
\end{figure}

In conclusion, our numerical results validate the quality of the lower and
upper bounds, $\beta _{2}^{\ast }$ and $\delta _{2}^{\ast }$, on the
spectral radius $\lambda _{1}$ in several social and communication networks.
Our bounds provide an interval of values in which the largest eigenvalue is 
\emph{guaranteed} to lie. This is in contrast with estimators based on
random networks, which can be very misleading and present no quality
guarantees.

\medskip

\section{Conclusions}

A fundamental question in the field of mathematical epidemiology is to
understand the relationship between a network's structural properties and
its epidemic threshold. For many virus epidemic models, the role of the
network topology is characterized by the largest eigenvalue of its adjacency
matrix, such that the larger the eigenvalue, the more efficient a network is
to spread a disease (or a piece of information) virally. In many cases of
practical interest, it is not possible to retrieve the complete structure of
a network of contacts due to privacy and/or security constrains. Thus, it is
not possible to exactly compute the largest eigenvalue of the network. On
the other hand, it is usually easy to retrieve local views of a network,
also called egonets, by extracting the structure of neighborhoods around a
collection of chosen nodes. To estimate the value of the spectral radius
when only egonets are available, researchers usually use random network
models in which they prescribe local structural features that can be
extracted from the egonets, such as the degree distribution. This approach,
although very common in practice, presents a major flaw: Random network
models implicitly induce many structural properties that are not directly
controlled and can be relevant to the spreading dynamics.

In this paper, we have presented an alternative mathematical framework,
based on algebraic graph theory and convex optimization, to study how
egonets constrain the interval of possible values in which the largest
eigenvalue (and, therefore, the epidemic threshold) must lie. Our approach
provides an interval of values in which the largest eigenvalue is \emph{%
guaranteed} to lie and is applicable to weighted networks. This is in
contrast with estimators based on random networks, which can be very
misleading and present no quality guarantees. Our numerical simulations have
shown that the resulting interval in which the largest eigenvalue must lie
is very narrow for several social and communication networks. This indicates
that, for an important collection of networks, the viral epidemic threshold
is strongly constrained by local structural properties of the network.

\appendix

\section{Proof of Corollary \protect\ref{Three Moments Bound}}

\textbf{Corollary \ref{Three Moments Bound} }Let $\mathcal{G}$ be a simple
graph with adjacency matrix $A_{\mathcal{G}}$. Denote by $n$, $e$, and $%
\Delta $ the number of nodes, edges, and triangles in $\mathcal{G}$,
respectively. Then,%
\begin{equation*}
\lambda _{1}\left( A_{\mathcal{G}}\right) \geq \frac{3\Delta +\sqrt{9\Delta
^{2}+8e^{3}/n}}{2e}.
\end{equation*}

\begin{proof}
From Corollary \ref{Three Moments}, we have that the first three moments of $%
\mathcal{G}$ are $m_{1}(A_{\mathcal{G}})=0,~m_{2}(A_{\mathcal{G}})=2e/n,~$%
and $m_{3}(A_{\mathcal{G}})=6\Delta /n$ (by definition, $m_{0}(A_{\mathcal{G}%
})=1$). Substituting the sequence of moments, $\mathbf{m}_{3}=\left(
1,0,2e/n,6\Delta /n\right) $, into (\ref{SDP generic bound}), we have that $%
\beta _{1}^{\ast }\left( \mathbf{m}_{3}\right) $ is the solution to the
following SDP:%
\begin{equation*}
\begin{array}{rl}
\min & x \\ 
\text{s.t.} & R\left( x\right) \triangleq \left[ 
\begin{array}{cc}
x & -2e/n \\ 
-2e/n & 2ex/n-6\Delta /n%
\end{array}%
\right] \succcurlyeq 0.%
\end{array}%
\end{equation*}%
The characteristic polynomial of $R\left( x\right) $ can be written as $\phi
\left( s;x\right) =\det \left( sI-R\left( x\right) \right) =s^{2}-s~tr\left(
R\left( x\right) \right) +\det \left( R\left( x\right) \right) $. Then, $%
R\left( x\right) \succcurlyeq 0$, if and only if both roots of $R\left(
x\right) $ are nonnegative. By Descartes' rule, this happens if and only if
the following two conditions are satisfied:

\begin{enumerate}
\item[(1)] $tr\left( R\left( x\right) \right) =x\left( 1+2e/n\right)
-6\Delta /n\geq 0$, which implies 
\begin{equation}
x\geq \frac{6\Delta }{2e+n}\triangleq x_{1}.  \label{Cond1}
\end{equation}

\item[(2)] $\det \left( R\left( x\right) \right) =2ex^{2}/n-6\Delta
x/n-4e^{2}/n^{2}\geq 0$, which implies%
\begin{equation}
x\geq \frac{3\Delta +\sqrt{9\Delta ^{2}+8e^{3}/n}}{2e}\triangleq x_{2}.
\label{Cond2}
\end{equation}
\end{enumerate}

We also have that, $x_{2}>\frac{3\Delta +\sqrt{9\Delta ^{2}}}{2e}=\frac{%
3\Delta }{e}>\frac{6\Delta }{2e+n}=x_{1}.$Therefore, the minimum value of $x$
satisfying (\ref{Cond1}) and (\ref{Cond2}) is equal to the right hand side
of (\ref{Bound from three}).
\end{proof}


\end{document}